\documentclass[12pt]{article}
\begin{document}

\begin{titlepage}

\vspace*{-2cm}

\vspace{.5cm}

\begin{centering}

\huge{$L_{\infty}$ algebra structures of Lie algebra deformations }

\vspace{.5cm}

\large  {Jining Gao }\\

\vspace{.5cm}

Department of Mathematics, North Carolina State University,
Raleigh, NC 27695-8205.

\vspace{.5cm}

\begin{abstract}

In this paper,we will show how to kill the obstructions to Lie
algebra deformations via a method which  essentially embeds a Lie algebra into Strong homotopy
Lie algebra or $L_{\infty}$ algebra. All such obstructions have been transfered to the revelvant
$L_{\infty}$ algebras which contain only three terms.  \end{abstract}

\end{centering}

\end{titlepage}

\pagebreak

\def\lh{\hbox to 15pt{\vbox{\vskip 6pt\hrule width 6.5pt height 1pt}
  \kern -4.0pt\vrule height 8pt width 1pt\hfil}}
\def\blob{\mbox{$\;\Box$}}
\def\qed{\hbox{${\vcenter{\vbox{\hrule height 0.4pt\hbox{\vrule width
0.4pt height 6pt \kern5pt\vrule width 0.4pt}\hrule height 0.4pt}}}$}}

\newtheorem{theorem}{Theorem}
\newtheorem{lemma}[theorem]{Lemma}
\newtheorem{definition}[theorem]{Definition}
\newtheorem{corollary}[theorem]{Corollary}
\newtheorem{proposition}[theorem]{Proposition}
\newcommand{\proof}{\bf Proof.\rm}

\section{Introduction}

In the last two decades,  deformations of various types of structures
have assumed an ever
increasing role in mathematics and physics. For each such deformation
problem  a goal  is to determine
if all related deformation obstructions  vanish and many beautiful
techniques been
developed  to determine when this is so. Sometimes genuine deformation
obstructions arise and occasionally  that closes mathematical
development in such cases, but
in physics such problems are dealt with by
introducing  new auxiliary fields to kill such obstructions. This idea
suggests that one might deal with deformation problems  by enlarging the
relevant category to a new category obtained by
appending additional algebraic structures to the old category.
To achieve the purpose of removing obstructions to Lie algebra deformation,
we embed
the Lie algebra
into an appropriate sh-Lie algebra  in such a way that the obstructions
will vanish
in the category of sh-Lie algebra deformations.
  In order to be complete we review  basic facts on Lie
algebra deformations; more detail may be found in the book edited by
M.Hazawinkel
and M.Gerstenhaber \cite{HG}.

\section{Deformation theory, sh-Lie algebras }

Let $A$ be a $k$-algebra and $\alpha$ be its multiplication, i.e.,
$\alpha$ is  a
$k$-bilinear map $A\times A \longrightarrow A$ defined by
$\alpha(a,b)=ab.$ A deformation of $A$ may be defined to be a formal
power series
$\alpha_t=\alpha+t\alpha_1+t^2\alpha_2+\cdots $
where each $\alpha_i:A\times A \longrightarrow A$ is a $k$-bilinear map
and the
``multiplication" $\alpha_t$ is formally of
the same "kind" as $\alpha,$ e.g., it is  associative or Lie or whatever
is required. One
technique used  to set up a deformation problem is to extend a
k-bilinear mapping
$\alpha_t:A\times A \longrightarrow A[[t]]$ to  a  $k[[t]]$-bilinear
mapping $\alpha_t:A[[t]]\times A[[t]]
\longrightarrow A[[t]].$ A mapping
$\alpha_t:A[[t]]\times A[[t]] \longrightarrow A[[t]]$ obtained in this
manner  is
necessarily uniquely determined by it's
values on $A\times A.$ In fact we would not regard the mapping
$\alpha_t:A[[t]]\times A[[t]] \longrightarrow A[[t]]$
to be a deformation of $A$ unless it is determined by it's values on $
A\times A$.

From this point on, we assume that $(A,\alpha)$ is a Lie algebra,i.e.,
we assume that
$\alpha(\alpha(a,b),c)+\alpha(\alpha(b,c),a)+\alpha(\alpha(c,a),b)=0$.
Thus the problem of
deforming a Lie algebra $A$ is equivalent to the problem of finding a
mapping
$\alpha_t:A\times A \longrightarrow A[[t]]$ such that
$\alpha_t(\alpha_t(a,b),c)+\alpha_t(\alpha_t(b,c),a)+\alpha_t(\alpha_t(c,a),b)=0.$
If we set $\alpha_0=\alpha$ and expand this Jacobi identity by making
the substitution
$\alpha_t=\alpha+t\alpha_1+t^2\alpha_2+\cdots ,$ we get the equation
\begin{eqnarray}
\sum_{i,j=0}^{\infty}[\alpha_j(\alpha_i(a,b),c)+\alpha_j(\alpha_i
(b,c),a)+\alpha_t(\alpha_i (c,a),b)]t^{i+j}=0
\end{eqnarray}
and consequently a sequence of deformation equations;
\begin{eqnarray}
\sum_{i,j\geq 0,i+j=n}[\alpha_j(\alpha_i(a,b),c)+\alpha_j(\alpha_i
(b,c),a)+\alpha_t(\alpha_i (c,a),b)]=0.
\end{eqnarray}
The first two equations are:
\begin{eqnarray}
\alpha_0(\alpha_0(a,b),c)+\alpha_0(\alpha_0(b,c),a)+\alpha_0(\alpha_0(c,a),b)=0
\\
\alpha_0(\alpha_1(a,b),c)+\alpha_0(\alpha_1(b,c),a)+\alpha_0(\alpha_1(c,a),b)+\alpha_1(\alpha_0(a,b),c)\nonumber
\\
  +\alpha_1(\alpha_0(b,c),a)+\alpha_1(\alpha_0(c,a),b)=0
\end{eqnarray}

We can reformulate the discussion above in a slightly more compact form.
Given
a sequence
$\alpha_n: A\times A\longrightarrow A$
of bilinear maps, we define ``compositions" of various of the
$\alpha_n$ as follows:
\begin{eqnarray}
& & \alpha_i\alpha_j:A\times A \times A \longrightarrow A
\end{eqnarray}
is defined by
\begin{eqnarray}
& & (\alpha_i \alpha_j)(x_1,x_2,x_3)=\sum_{\sigma \in
unsh(2,1)}(-1)^{\sigma} \alpha_i (\alpha_j(x_{\sigma(1)},
x_{\sigma(2)}),x_{\sigma(3)})
\end{eqnarray}
for arbitrary $x_1,x_2,x_3\in A.$

Thus the deformation equations are equivalent to following equations:
\begin{eqnarray}
\alpha_{0}^{2}=0 \\ \alpha_0\alpha_1+\alpha_1\alpha_0=0 \\
\alpha_1^2+\alpha_0\alpha_2+\alpha_2\alpha_0=0\\
\cdots \nonumber\\ \Sigma_{i+j=n} \alpha_i\alpha_j=0\\ \cdots.
\end{eqnarray}

Define a bracket on the sequence $\{\alpha_n\}$ of mappings by
$[\alpha_i,\alpha_j]=\alpha_i\alpha_j+\alpha_j\alpha_i$ and a
``differential"
$d$ by $d=ad_{\alpha_0}=[\alpha_0, \cdot ],$ the ``adjoint
representation" relative
to $\alpha_0.$  Notice that the second equation in the list above is
equivalent to the
statement that $\alpha_1$ defines a cocycle $\alpha_1\in Z^2(A,A)$ in
the Lie
algebra cohomology of $A.$
Moreover it is known that the second cohomology group $H^2(A,A)$
classifies the equivalence class of infinitesimal deformations of $A$
\cite{HG}.
This being the case we refer to the triple $(A,\alpha_0,\alpha_1)$ as
being initial
conditions for deforming the Lie algebra $A.$
Notice that the third equation in the above list can be rewritten as
\begin{eqnarray}
[\alpha_1,\alpha_1]=-[\alpha_0,\alpha_2]=-d\alpha_2
\end{eqnarray}
When this  equation holds one has then that  $[\alpha_1,\alpha_1]$ is a
coboundary
and so defines
the trivial element of $H^3(A,A)$ for any given deformation $\alpha_t.$
Thus if
$[\alpha_1,\alpha_1]$ is not a coboundary,  then we may regard
$[\alpha_1,\alpha_1]$ as the first  obstruction to deformation and in
this case we
can not deform
$A$ at second order. In general, to say that there exists a  deformation of
$(A,\alpha_0,\alpha_1)$ up to order $n-1,$ means that  there
exists a sequence of maps $\alpha_0,\cdots, \alpha_{n-1}$ such that
$\sum_{\sigma \in unsh(2,1)}(-1)^{\sigma} \alpha_t  (\alpha_t(x_{\sigma(1)},
x_{\sigma(2)}),x_{\sigma(3)})=0  \quad (mod \quad t^n). $
If this is the case and if there is an obstruction to  deformation  at
$n$th order,
then it follows  that $\rho_n=-(\Sigma_{i+j=n,i,j>0}
\alpha_i\alpha_j)$ is in some sense the obstruction and $[\rho_n]$ is a
nontrivial
element of $H^3(A,A).$  If $[\rho_n]\neq 0$,then the process of obtaining a
deformation will terminate at order $n-1$ due to the existence of the
obstruction
$\rho_n.$  In principal, it is possible that one could return to the
beginning and
select different
terms for the $\alpha_i$ but when this fails what can one say? This is
the issue in
the remainder
of this section.

Indeed  the central point of this section  is to show
that when there is an obstruction to the deformation of a Lie algebra,
one can
use
the obstruction itself to define one of the structure mappings of an
sh-Lie algebra.
Without loss of generality, we  consider a  deformation problem which
has a first
order obstruction.

The required sh-Lie structure lives on a graded vector space $X_*$ which
we define
below. This space in degree zero is given by $X_0=A[[t]]=(\{\Sigma
a_it^i |\quad a_i\in A\}.$ The spaces  $\cal B=$
$ <t^2>=$ $A[[t]]\cdot t^2=$ $\{\Sigma_{i\geq 2} a_it^i |a_i\in A\}$ and
$\cal F= $
$X_0/{\cal B}=$  are also relevant to our construction. Notice that
$\cal F$ is isomorphic to $\{a_0+a_1 t | \quad a_0,a_1 \in A\}$ as a
linear space and that   $X_0,\cal B$ are both $k[[t]]$-modules while
$\cal F$ is a $k[[t]]/<t^2>$ module (recall that $k$ is
underlying field of $A$).
To summarize, we have
following short exact sequence:
$$
0{\longrightarrow} {\cal B} {\longrightarrow} {X}_{0} {\longrightarrow
}{\cal F}
\longrightarrow{0}.
$$
Suppose that the initial Lie structure of $A$ is given  by
$\alpha_0: A\times A \longrightarrow A$ and denote a fixed
infinitesimal deformation  by $[\alpha_{1}] \in H^{2}(A,A).$  One of the
structure mappings
of our sh-Lie structure will be determined by the mapping $\tilde l_2 :
X_0 \times X_{0} \longrightarrow X_{0} $
defined as follows: for any $a,b \in A$, let
\begin{eqnarray}
\tilde l_2(a,b)=\alpha_0(a,b)+\alpha_1(a,b)t
\end{eqnarray}
and extend it to $X_0$ by requiring that it be  $k[[t]]$-bilinear.
Obviously, $\tilde
l_2$ induces a Lie bracket $[,]$ on $\cal F$, but if the  obstruction
$[\alpha_1,\alpha_1] $ is not zero,
then $\tilde l_2^2 \neq 0$ and consequently $\tilde l_2$ can not be a
Lie bracket on $X_0$
(since it doesn't satisfy the Jacobi identity).

To deal with this obstruction we will show that we can use
$\alpha_0,\alpha_1$ to
construct an
sh-Lie structure with at most three nontrivial structure maps
$l_1,l_2,l_3$ such
that the value of $l_3$ on $A\times A\times A$
is the same as that of $[\alpha_1,\alpha_1].$ In
particular, $l_3$ will vanish if and only if the obstruction
$[\alpha_1,\alpha_1]$
vanishes. Thus the sh-Lie algebra encodes the obstruction to deformation
of the Lie
algebra $(A,\alpha_0).$

The required sh-Lie algebra lives on a certain homological resolution
$(X_*,l_1)$ of
${\cal F},$
so our first task is to construct this resolution space for $ {\cal F}.$ To
do this let's introduce a ``superpartner set of $A,$" denoted by  $A[1],$ as
follows:
for each $a\in A,$ introduce $a^*$ such that $a^*\leftrightarrow a$ is a
one to one
correspondence and define $\epsilon(a^*)=\epsilon(a)+1.$
Let $X_1=A[1][[t]]t^2$ and define a map $l_1: X_1 \longrightarrow X_0$ by
$$ l_1(x)=\Sigma_{i\geq 2} a_i t^i \in X_0, \quad \quad
x=\Sigma_{i\geq 2}a_i^*t^i \in X_1 .      $$

Notice that this is just the $k[[t]]$ extension of the
$a^*\leftrightarrow a$  map.
Since $l_1$ is injective, we obtain a homological resolution
$X_*=X_0\oplus X_1$
due to the fact that  the  complex defined by:
\begin{eqnarray}
0\longrightarrow X_1\stackrel{l_1}{\longrightarrow}X_0\longrightarrow 0
\end{eqnarray}
has the obvious property that $H(X_*)=H_0(X_*)\simeq \cal F.$

The $sh$-Lie algebra being constructed will have the property that
$l_n =0,n\geq 4.$ Generally sh-Lie algebras can have any number of
nontrivial
structure maps. The fact that all the structure mappings of our sh-Lie
algebra are
zero with the exception of $l_1,l_2,l_3$ is an immediate consequece of
the fact
that we are able to produce a
resolution of the space $\cal F$ such that $ X_k =0 $ for $k \geq 2.$ In
general
such resolutions do not exist and so one does not have $ l_n =0 $ for
$n\geq 4.$

In order to finish the
preliminaries, we now  construct a contracting homotopy $s$
such that following commutative diagram holds:
$$
\begin{array}{ccccc}
&&s&& \\
0\longrightarrow & X_1 &\stackrel
{\textstyle\longleftarrow}{\longrightarrow} &X_0&\longrightarrow 0\\
&&l_1&& \\ 
& \lambda\Bigl\uparrow\Bigr\downarrow\eta & &
\lambda\Bigl\uparrow\Bigr \downarrow\eta & \\
&&&& \\
0\longrightarrow &0&\longrightarrow & \cal F & \longrightarrow 0
\end{array}
$$
Clearly the  linear space $X_0$ is the direct sum of ${\cal B}$ and
a complementary subspace which is isomorphic to ${\cal F};$ consequently
we have
$X_0 \simeq {\cal B} \oplus \cal F.$  Define $\eta=proj\mid_{\cal F},
\lambda=i_{{\cal
F}\rightarrow X_0}$
and a contracting homotopy $s: X_0\longrightarrow X_1$ as follows:
write $X_0={\cal B}\oplus \cal F,$  set $s|_{\cal F}=0,$ and let $
s(x)=-x^*$ for
all $ x\in {\cal B}.$
It is easy to show that
$ \lambda \circ \eta  -1_{X_*}=l_1\circ s +s \circ l_1.$
In order to obtain the sh-Lie algebra referred to above, we apply a
theorem of
\cite{BFLS}. The hypothesis of this theorem requires the existence of a
bilinear mapping
$\tilde l_2$ from $X_0\times X_0$ to $X_0$ with the properties that  for
$c,c_1,c_2,c_3\in X_0$ and $b\in {\cal B}$
$(i) \quad  \tilde l_2(c,b)\in {\cal B}$ and $(ii) \quad  \tilde
l_2^2(c_1,c_2,c_3)\in {\cal B}.$
To see that $(i)$ holds notice that if $p(t),q(t)\in X_0=A[[t]],$ then
$\tilde l_2(p(t),q(t)t^2)=r(t)t^2$ for some $r(t)\in A[[t]]=X_0.$
Also note that the fact that $\tilde l_2$ induces a Lie bracket on ${\cal
F}=X_0/{\cal B}$ implies that
$\tilde l_2^2$ is zero modulo ${\cal B}$ and $(ii)$ follows. Thus $X_*$
supports an
sh-Lie structure
with only three nonzero structure maps $l_1,l_2,l_3$ (see the remark at
the end of
\cite{BFLS}).

\begin{theorem} Given a Lie algebra $A$ with Lie bracket $\alpha_0$ and an
infinitesimal obstruction $[\alpha_1]\in H^2(A,A)$ to deforming
$(A,\alpha_0),$ there
is an sh-Lie algebra on the graded space
$(X_*,l_1)$ with structure maps $\{l_i\}$ such that $l_n=0$ for $n\geq
4.$ The
graded space
$X_*$ has at most two nonzero terms $X_0=A[[t]],X_1=A[1][[t]]t^2.$
Finally, the maps
$l_1,l_2,l_3$ may be given explicitly in terms of the maps
$\alpha_0,\alpha_1.$
\end{theorem}

${\bf Remark:}$ The mapping $l_1$ is simply the differential of the
graded space $(X_*,l_1).$ The mapping $l_2$
restricted
to $X_0 \times X_0$ is the mapping $\tilde l_2$
defined directly in terms of $\alpha_0,\alpha_1$ above.  On $X_1\times
X_0,$ $l_2$ is
determined by
$l_2(a^*t^2,b) =t^2(\alpha_0(a,b)^*+\alpha_1(a,b)^*t)$ for $a^*\in
A[1],b\in A.$
Finally, $l_3$ is uniquely determined by its values on $A\times A\times
A\subset
X_0\times X_0\times X_0$ and is explicitly a multiple of the obstruction
to the
deformation
of $(A,\alpha_0), $ in particular,
$  l_3(a_1,a_2,a_3)=-t^2[\alpha_1,\alpha_1](a_1,a_2,a_3),  a_i\in A.$

\proof { The sh-Lie structure maps are given by Theorem 7 of
\cite{BFLS}. The fact that
$l_n=0,n\geq 4$ is an observation of Markl which was proved by Barnich
\cite{B}
( see the remark at the end of \cite{BFLS}).
A generalization of Markl's remark is available in a paper by Al-Ashhab
\cite{Samer}
and in that paper more explicit formulas are given for $l_1,l_2,l_3.$
Examination of these
formulas  provide the details needed for the calculations below.

First of all, we examine the mapping
$l_2: X_*\times X_*\longrightarrow X_*.$ Now
$ l_2 : A \times A $ $\longrightarrow X_*$
is  determined by $\tilde l_2: X_0\times X_0\longrightarrow X_0,$
consequently we
need only consider the restricted mapping:
\begin{eqnarray}
l_2: X_1\times X_0\longrightarrow X_1.
\end{eqnarray}
Moreover, since $X_0$ is a module over $k[[t]],$  $X_1$ is a module over
$k[[t]]t^2,$
and $\tilde l_2$ respects these structures
we need only consider its values on pairs $(a^{*}t^2,b)$ with
$a^{*}t^2\in X_1,b\in X_0.$   By  Theorem 2.2 of \cite{Samer}, we have
\begin{eqnarray}
l_2(a^{*}t^2,b)=-sl_2 l_1[(a^{*}t^2)\otimes b]\nonumber\\
=-sl_2[l_1(a^{*}t^2)\otimes b+(-1)^{\epsilon(a^*)}(a^{*}t^2)\otimes
l_1(b)]\nonumber\\
=-sl_2[(at^{2}\otimes b)]=-s[t^{2} l_{2}(a \otimes b)] \nonumber \\
=-s[t^2(\alpha_0(a,b)+\alpha_1(a,b)t)]\nonumber\\
=-s[\alpha_0(a,b)t^2+\alpha_1(a,b)t^3]\nonumber\\
=\alpha_0(a,b)^*t^2+\alpha_1(a,b)^*t^3\nonumber\\
=t^2(\alpha_0(a,b)^*+\alpha_1(a,b)^*t).
\end{eqnarray}

From this deduction, we  that the mapping $ l_{2} $ can essentially be
replaced  by
the modified map:
\begin{eqnarray}
\bar l_2: A[1] \times A\longrightarrow A[1][[t]],  \quad \quad \bar
l_2(a^*,b)=\alpha_{0}(a,b)^*+\alpha_{1}(a,b)^*t.
\end{eqnarray}
We clarify this remark below by showing that a new sh-Lie structure can
be obtained
with
$\bar l_2$ playing the role of $l_2.$

The next  mapping we examine is the mapping
\begin{eqnarray}
l_3: X_0 \times X_0\times X_0 \longrightarrow X_1
\end{eqnarray}
Since $l_{3} $ is $k[[t]]$-linear,we need only consider  mappings of the
type:
\newline $l_3: A\times A\times A \longrightarrow X_1$
where for $ x_1,x_2,x_3\in A,$
\begin{eqnarray}
l_3(x_1,x_2,x_3)=sl_2^2(x_1,x_2,x_3)\nonumber\\
=\sum_{\sigma \in unsh(2,1)}(-1)^{\sigma} sl_2(l_2(x_{\sigma(1)},
x_{\sigma(2)}),x_{\sigma(3)})\nonumber\\
=\sum_{\sigma \in unsh(2,1)}(-1)^{\sigma} sl_2
(\alpha_0(x_{\sigma(1)}, x_{\sigma(2)})+\alpha_1(x_{\sigma(1)},
x_{\sigma(2)})t,x_{\sigma(3)})\nonumber\\
  =\sum_{\sigma \in unsh(2,1)}(-1)^{\sigma}
s[\alpha_0(\alpha_0(x_{\sigma(1)},
x_{\sigma(2)}),x_{\sigma(3)})+\nonumber\\
   t\alpha_1(\alpha_0(x_{\sigma(1)}, x_{\sigma(2)}),x_{\sigma(3)})+
t\alpha_0(\alpha_1(x_{\sigma(1)}, x_{\sigma(2)}),x_{\sigma(3)})\nonumber\\
   +t^2 \alpha_1(\alpha_1(x_{\sigma(1)},
x_{\sigma(2)}),x_{\sigma(3)}),x_{\sigma(3)})\nonumber\\
  =s(\sum_{\sigma \in unsh(2,1)}(-1)^{\sigma}
\alpha_0(\alpha_0(x_{\sigma(1)},
x_{\sigma(2)}),x_{\sigma(3)})+\nonumber\\
+ t(\sum_{\sigma \in unsh(2,1)}(-1)^{\sigma}
\alpha_1(\alpha_0(x_{\sigma(1)},
x_{\sigma(2)}),x_{\sigma(3)})+\nonumber\\
\sum_{\sigma \in unsh(2,1)}(-1)^{\sigma} \alpha_0(\alpha_1(x_{\sigma(1)},
x_{\sigma(2)}),x_{\sigma(3)}))\nonumber\\
+t^2(\sum_{\sigma \in unsh(2,1)}(-1)^{\sigma}
\alpha_1(\alpha_1(x_{\sigma(1)},
x_{\sigma(2)}),x_{\sigma(3)}))]\nonumber\\
  =s((\alpha_0^2+t(\alpha_0\alpha_1+\alpha_1\alpha_0)+t^2
\alpha_1^2)(x_1,x_2,x_3))\nonumber\\
=s(t^2 \alpha_1^2 (x_1,x_2,x_3))\nonumber\\
  =-t^2 (\alpha_1^2 (x_1,x_2,x_3))^*
\end{eqnarray}
Or $l_3(x_1,x_2,x_3)=-t^2 ([\alpha_1,\alpha_1] (x_1,x_2,x_3))^*$ which
is precisely the
``first deformation obstruction class".}

Recall that we  know from Theorem 7 of \cite{BFLS} that we have an
sh-Lie structure.
The point of these calculations is that it enables us to obtain the
modified sh-Lie structure of Corollary 10 below and it is this structure
which is relevant to Lie algebra deformation.
Thus we already know that the mappings $l_1,l_2,l_3$  satisfy the relations:
\begin{eqnarray}
l_1l_2-l_1l_2=0   \\ {l_2}^2+l_1l_3+l_3l_1=0   \\{l_3}^2=0  \\
l_2l_3+l_3l_2=0.
\end{eqnarray}

Observe that if we let $\tilde X_*=\tilde X_1 \oplus \tilde X_0=A[1] [
[t] ] \oplus
A[[t]],$ then the formulas
defining $l_1,l_2, l_3$  defined on $X_*$ make sense on the new complex
$\tilde
X_*.$ Indeed
the calculations above show that  $l_1,l_3$ are uniquely determined by
their values
on ``constants" in the sense that they could be first defined on elements of
$A[1]\oplus A \subseteq A[1][[t]]\oplus A[[t]]$ and then extended to
$A[1][[t]]\oplus
A[[t]]$ using the fact that $l_1,l_3$ are required to be $k[[t]]$
linear. $l_2$ is not
obviously $k[[t]]$ linear. The whole point of corollary 10 below is that
the $sh$-Lie
structure defined by Theorem 10 can be redefined to obtain $sh$-Lie maps on
the graded space $\tilde X_*$ which are obviously
$k[[t]]$ linear and consequently this "new" structure is intimately
related to
deformation theory.
Thus, as we say above, the modified map
$\bar l_2$ can be extended to the new complex $\tilde X_*$ and is uniquely
determined by its values on ``constants".
If we denote the extensions of $l_1,l_3$ to $\tilde X_*$ by $\bar
l_1,\bar l_3,$
then clearly these
mappings satisfy the same relations (63)-(66) as the maps $l_1,l_2,l_3$ and
consequently
if we define $\bar l_n=0, n\geq 4$ it follows that  $(\tilde X_*,\tilde
l_1,\tilde
l_2,\tilde l_3,0,0 \cdots) $ is an $sh$-Lie algebra. This proves the
following
corollary.

\begin{corollary} There is an sh-Lie structure on $A[1][[t]]\oplus
A[[t]]$ whose
structure mappings
$\{\bar l_1,\bar l_2,\bar l_3, 0,\cdots\} $are precisely the mappings
$\{l_1,l_2,l_3, 0,\cdots\} $ when restricted to  $A[1][[t]]t^2\oplus
A[[t]].$ Moreover,
the structure mappings of $A[1][[t]]\oplus A[[t]]$ have the property
that they are
uniquely determined by their values on $A[1]\oplus A$ and $k[[t]]$
linearity.
\end{corollary}

From the  discussion above
the set of mappings  $\{\bar  l_1,\bar l_2,\bar l_3\}$ is essentially a
deformation
of an $sh$-Lie algebra. In addition, the construction    of
the mapping  $\bar l_2$ is  equivalent to defining an initial condition
for a  Lie
algebra deformation.

This means that a Lie algebra which can't be deformed in the category of
Lie algebra
may admit an $sh$-Lie algebra deformation by first imbedding it into an
appropriate
$sh$-Lie algebra.

\section*{Acknowledgments}
I would like to thank R.Fulp for useful  discussions throughout the writing
of this paper with special thanks
for his revision of the manuscript. I am also pleased to express my 
appreciation to
J.Stasheff who provided many thoughtful comments which led to substantial
improvement of  an early draft of my paper.

\newpage

\ifx\undefined\bysame
\newcommand{\bysame}{\leavevmode\hbox to3em{\hrulefill}\,}
\fi

\end{document}